\documentclass{amsart}

\DeclareMathSymbol{\twoheadrightarrow}
{\mathrel}{AMSa}{"10}

            \def\H{{\mathbf H}}
\def\Q{{\mathbf Q}}
\def\Z{{\mathbf Z}}

\def\RR{{\mathfrak R}}
\def\F{{\mathbf F}}
\def\SS{{\mathbf S}}
\def\A{{\mathbf A}}
\def\Sn{{\mathbf S}_n}
\def\An{{\mathbf A}_n}

\def\Gal{\mathrm{Gal}}
\def\Perm{\mathrm{Perm}}

\def\End{\mathrm{End}}
\def\Aut{\mathrm{Aut}}
\def\Hom{\mathrm{Hom}}

\def\I{\mathrm{Id}}

\def\fchar{\mathrm{char}}

\def\GL{\mathrm{GL}}

\def\L{\mathrm{L}}

\def\SL{\mathrm{SL}}
\def\PSL{\mathrm{PSL}}

\def\Sp{\mathrm{Sp}}
\def\M{\mathrm{M}}
\def\dim{\mathrm{dim}}

\def\P{{\mathbf P}}

\def\a{{\mathfrak a}}

\newtheorem{thm}{Theorem}[section]
\newtheorem{lem}[thm]{Lemma}

\theoremstyle{definition}

\newtheorem{ex}[thm]{Example}

\newtheorem{rem}[thm]{Remark}

\title[Non-supersingular Hyperelliptic jacobians]
{Non-supersingular Hyperelliptic jacobians}
\author[Yuri G. Zarhin]{Yuri G. Zarhin}
\address{Department of Mathematics, Pennsylvania State University,
University Park, PA 16802, USA} \email{zarhin\char`\@math.psu.edu}

\begin{document}
\begin{abstract}
In his previous papers ~\cite{ZarhinMRL,ZarhinMRL2,ZarhinMMJ} the
author proved that in characteristic $\ne 2$ the jacobian $J(C)$
of a hyperelliptic curve $C: y^2=f(x)$ has only trivial
endomorphisms over an algebraic closure $K_a$ of the ground field
$K$ if the Galois group $\Gal(f)$ of the irreducible polynomial
$f(x) \in K[x]$ is either the symmetric group $\Sn$ or the
alternating group $\A_n$. Here $n\ge 9$ is the degree of $f$. The
goal of this paper is to extend this result to the case when
either $n=7,8$ or $ n=5,6$ and $\fchar(K)>3$.


%11G30;11G10.

\end{abstract}
\maketitle
\section{Introduction}
Let $K$ be a field and $K_a$ its algebraic closure. Assuming that
$\fchar(K)=0$, the author \cite{ZarhinMRL} proved that  the
jacobian $J(C)=J(C_f)$  of a hyperelliptic curve
$$C=C_f:y^2=f(x)$$
has only trivial endomorphisms over  $K_a$ if the Galois group
$\Gal(f)$ of the irreducible polynomial $f \in K[x]$ is ``very
big". Namely, if $n=\deg(f) \ge 5$ and $\Gal(f)$ is either the
symmetric group $\Sn$ or the alternating group $\An$
 then the ring $\End(J(C_f))$ of $K_a$-endomorphisms of $J(C_f)$ coincides with $\Z$.
 Later the author ~\cite{ZarhinMRL,ZarhinMMJ} extended
 this result to the case of positive $\fchar(K)>2$ but under the additional assumption
that $n \ge 9$, i.e., the genus of $C_f$ is greater or equal than
$4$. We refer the reader to
~\cite{Mori1,Mori2,Katz1,Katz2,Masser,KS,ZarhinMRL,ZarhinMRL2,ZarhinMMJ,ZarhinP,{ZarhinPAMS}}
for a discussion of known results about, and examples of,
hyperelliptic jacobians without complex multiplication.

The aim of the present paper is to extend this result to the case
when either $n \ge 7$ or when $n \ge 5$ but $\fchar(K)>3$. Notice
that it is known ~\cite{ZarhinMRL} that in those cases either
$\End(J(C))=\Z$ or $J(C)$ is a supersingular abelian variety and
the real problem is how to prove that $J(C)$ is {\sl not}
supersingular.

We also discuss  the case of two-dimensional
 $J(C)$ in characteristic $3$.

\section{Main result}
\label{mainr} Throughout this paper we assume that $K$ is a field
of  characteristic $p$ different from $2$. We fix its algebraic
closure $K_a$ and write $\Gal(K)$ for the absolute Galois group
$\Aut(K_a/K)$.

\begin{thm}
\label{main}
Let $K$ be a field with $p=\fchar(K) > 2$,
 $K_a$ its algebraic closure,
$f(x) \in K[x]$ an irreducible separable polynomial of  degree
$n$. Let us assume  that $\Gal(f)=\Sn$ or $\An$. Suppose that $n$
enjoys one of the following properties:
\begin{enumerate}
\item[(i)]
$n=7$ or $8$;
\item[(ii)]
$n=5$ or $6$. In addition, $p=\fchar(K)>3$.
\end{enumerate}

 Let $C_f$ be the hyperelliptic curve
$y^2=f(x)$. Let $J(C_f)$ be its jacobian, $\End(J(C_f))$ the ring
of $K_a$-endomorphisms of $J(C_f)$. Then $\End(J(C_f))=\Z$.
\end{thm}

\begin{rem}
\label{redA}
 Replacing $K$ by a suitable finite separable
extension, we may assume in the course of the proof of Theorem
\ref{main} that $\Gal(f)=\An$. Taking into account that $\An$ is
simple non-abelian and replacing $K$ by its abelian extension
obtained by adjoining to $K$ all $2$-power roots of unity, we may
also assume that $K$ contains all $2$-power roots of unity.
\end{rem}

\begin{rem}
\label{odd} Let $f(x) \in K[x]$ be an irreducible separable
polynomial of  {\sl even} degree $n=2m \ge 6$ such that
$\Gal(f)=\Sn$. Let $\alpha \in K_a$ be a root of $f$ and
$K_1=K(\alpha)$ be the corresponding subfield of $K_a$. We have
$$f(x)=(x-\alpha) f_1(x)$$ with $f_1(x) \in K_1[x]$. Clearly,
$f_1(x)$
 is an irreducible separable
polynomial over $K_1$ of  degree $n-1=2m-1$, whose Galois group is
$\SS_{n-1}$. It is also clear that the polynomials
$$h(x)=f_1(x+\alpha), h_1(x)=x^{n-1}h(1/x) \in K_1[x]$$ are
irreducible separable of  degree $n-1$ with the same Galois group
$\SS_{n-1}$.

The standard substitution
$$x_1=1/(x-\alpha), y_1=y/(x-\alpha)^m$$
establishes a birational isomorphism between $C_f$ and a
hyperelliptic curve
$$C_{h_1}: y_1^2=h_1(x_1).$$
\end{rem}

In light of results of ~\cite{ZarhinTexel,ZarhinPAMS} and Remarks
\ref{redA} and \ref{odd}, our
 Theorem \ref{main} is an immediate corollary of the following auxiliary
statement.

\begin{thm}
\label{main2} Let $K$ be a field with $p=\fchar(K) > 2$,
 $K_a$ its algebraic closure,
$f(x) \in K[x]$ an irreducible separable polynomial of  degree
$n$. Let us assume  that $n$ and the Galois group $\Gal(f)$ of $f$
enjoy one of the following properties:
\begin{enumerate}
\item[(i)]
$n=5$  and $\Gal(f)=\A_5$;
\item[(ii)]
$n=7$  and $\Gal(f)=\A_7$. In addition, $p=\fchar(K)
> 3$;
\end{enumerate}

Let $C$ be the hyperelliptic curve $y^2=f(x)$  and let
 $J(C)$  be the jacobian of $C$.

Then $J(C)$ is not a supersingular abelian variety.
\end{thm}

We will prove Theorem \ref{main2} in Section \ref{pmain2}.

Throughout the paper we write $\End^0(X)$ for the endomorphism
algebra $\End(X)\otimes\Q$ of an abelian variety $X$ over an
algebraically closed field $F_a$. Recall \cite{ZarhinMRL} that the
semisimple $\Q$-algebra $\End^0(X)$ has dimension $(2\dim(X))^2$
if and only if $p:=\fchar(F_a)\ne 0$  and $X$ is a supersingular
abelian variety. We write $\H_p$ is the quaternion $\Q$-algebra
unramified exactly at $p$ and $\infty$. It is well-known that if
$X$  is a supersingular abelian variety in characteristic $p$ then
$\End^0(X)$ is isomorphic to the matrix algebra $\M_g(\H_p)$ of
size $g:=\dim(X)$ over $\H_p$. We will use freely these facts
throughout the paper.

\section{Proof of Theorem \ref{main2}}
\label{pmain2}

We deduce Theorem \ref{main2} from the following statement.

\begin{thm}
\label{main3} Let $K$ be a field with $p=\fchar(K) > 2$,
 $K_a$ its algebraic closure, Let $n=q$ be an odd prime,
$f(x) \in K[x]$ an irreducible separable polynomial of  degree
$q$. Let us assume  that  the Galois group $\Gal(f)$ of $f$ is
$\L_2(q):= \PSL_2(\F_q)$, and that it acts  doubly transitively on the roots
of $f$. Suppose that either $q=5$ or $q=7$. Let $C$ be the
hyperelliptic curve $y^2=f(x)$  and let
 $J(C)$  be the jacobian of $C$.

If $J(C)$ is  a supersingular abelian variety then $n=5$ and $p=3$.
\end{thm}

\begin{proof}[Proof of Theorem \ref{main2} (modulo Theorem \ref{main3})]
If $n=5$ then $\A_5 \cong \L_2(5)$ and we are done.

Suppose that $n=7$. It is well-known that the simple non-abelian
group $\L_2(7)\cong \L_3(2):=\PSL_3(\F_2)$ acts doubly
transitively on the $7$-element projective plane $\P^2(\F_2)$  and
therefore is isomorphic to a doubly transitive subgroup of $\A_7$.
Hence there exists a finite algebraic extension $K_1$ of $K$ such
that the Galois group of $f$ over $K_1$ is $\L_2(7)$ acting doubly
transitively on the roots of $f(x)$. Applying Theorem \ref{main3}
to $K_1$ and $f$, we conclude that if $3 \ne
\fchar(K_1)=\fchar(K)=p$ then $J(C)$ is not supersingular.
\end{proof}

The following results will be used in order to prove Theorem \ref{main3}.

\begin{lem}
\label{double}
Let $K$ be a field with $\fchar(K) \ne 2$
 $K_a$ its algebraic closure, $\Gal(K)=\Aut(K_a)$ the Galois group of $K$.
Let $f(x) \in K[x]$ be an irreducible separable polynomial of  odd degree
$n$. Let us assume  that  $n \ge 5$ and the Galois group $\Gal(f)$ of $f$ acts
doubly transitively on the roots of $f(x)$.
Let $C$ be the hyperelliptic curve $y^2=f(x)$  and let
 $J(C)$  be the jacobian of $C$. Let $J(C)_2$ be the group of points of order $2$ in $J(C)(K_a)$ viewed as $\F_2$-vector space provided with a natural structure of  $\Gal(K)$-module.
Then the image of $\Gal(K)$ in $\Aut_{\F_2}(J(C)_2)$ is isomorphic to $\Gal(f)$ and
$$\End_{\Gal(K)}(J(C)_2)=\End_{\Gal(f)}(J(C)_2)=\F_2.$$
\end{lem}

\begin{thm}
\label{supernot} Let $F$ be a field with characteristic  $p>2$
and assume that $F$ contains all $2$-power roots of unity. Let
$F_a$ be an algebraic closure of $F$.

Let $G \ne \{1\}$ be a finite  perfect group.

Suppose that $g$ is a positive integer, $X$ is a supersingular
$g$-dimensional abelian variety defined over $F$. Let $\End(X)$ is
the ring of all $F_a$-endomorphisms of $X$ and
$\End^0(X)=\End(X)\otimes\Q$.

Let us assume  that the image of $\Gal(F)$ in $\Aut(X_2)$ is
isomorphic to $G$ and the corresponding faithful representation
$$\bar{\rho}: G \hookrightarrow \Aut(X_2) \cong \GL(2g,\F_2)$$
satisfies
$$\End_G{X_2}=\F_2.$$

 Then there exists a surjective group
homomorphism
$$\pi_1: G_1 \twoheadrightarrow G$$
enjoying the following properties:
\begin{enumerate}
\item[(a)]
The group $G_1$ is a perfect finite group. The kernel of $\pi_1$
is an elementary abelian $2$-group. \item[(b)]
 One may lift $\bar{\rho}\pi_1: G_1\to  \Aut(X_2)$ to a faithful absolutely irreducible symplectic
  representation
$$\rho: G_1 \hookrightarrow \Aut_{\Q_2}(V_2(X))$$
  of $G_1$ over $\Q_2$ in such a way that the following conditions hold:
\begin{itemize}
\item
 The character $\chi$ of $\rho$ takes values in $\Q$;
\item
  $\rho(G_1)\subset (\End^0(X))^*$;
\item
 The homomorphism from the group algebra $\Q[G_1]$ to
$\End^0(X)$ induced by $\rho$ is surjective and identifies
$\End^0(X)\cong \M_g(\H_p)$ with the  direct summand of $\Q[G_1]$ attached to $\chi$.
\end{itemize}
\item[(c)]
$p$ divides the order of $G$ and $p \le 2g+1$.
\item[(d)]
Suppose that either every homomorphism from $G$ to $\GL(g-1,\F_2)$
is trivial or  the $G$-module $X_2$ is very simple in the sense of
\cite{ZarhinTexel,ZarhinMMJ,ZarhinVery}. Then $\ker{\pi_1}$ is a
central cyclic group of order $1$ or $2$.
\end{enumerate}
\end{thm}

\begin{lem}
\label{sl2} Let $p$ be an odd prime.  Let $q$ be an odd prime and
$\Gamma=\SL_2(\F_q)$ or $\PSL_2(\F_q)$. Suppose that $q=5$ or $7$
and let us put $g=\frac{q-1}{2}$. Suppose that $\Q[\Gamma]$
contains a direct summand isomorphic to the matrix algebra
$\M_{g}(\H_p)$. Then $p=3$ and $q=5$.
\end{lem}

Theorem \ref{supernot} and Lemmas will be proven in Sections
\ref{main3p} and \ref{pl}.

\begin{proof}[Proof of Theorem  \ref{main3} (modulo Theorem \ref{supernot} and Lemmas \ref{double} and \ref{sl2})]
Let us put $$X=J(C), G=\PSL_2(\F_q), g =\frac{q-1}{2}.$$ Clearly,
either $q=5,g=2$ or $q=7,g=3$. In both cases $g=\dim(X)$, the
group $G$ is simple and $\GL(g-1,\F_2)$ is solvable. It follows
that every homomorphism from $G$ to $\GL(g-1,\F_2)$ is trivial.
 It
follows from Lemma \ref{double} that the image of $\Gal(K)$ in
$\Aut(X_2)$ is isomorphic to $G$ and the corresponding faithful
representation
$$\bar{\rho}: G \hookrightarrow \Aut(X_2) \cong \GL(2g,\F_2)$$
satisfies
$$\End_G{X_2}=\F_2.$$
Let us assume that $X$ is supersingular. We need to get a
contradiction.

 Applying Theorem \ref{supernot}, we conclude that
there exist a finite perfect group $G_1$ and a surjective
homomorphism
$$\pi_1: G_1 \twoheadrightarrow G=\PSL_2(\F_q)$$
enjoying the following properties.
\begin{itemize}
\item[(i)]
Either $G_1 \cong G$ or $Z_1=\ker(\pi_1)$ is a central subgroup of
order $2$ in $G_1$;
\item[(ii)]
There exists a direct summand of $\Q[G_1]$ isomorphic to
$\M_g(\H_p))$.
\end{itemize}

The well-known description of central extensions of $\PSL_2(\F_q)$
when  $q$ is an odd prime \cite[\S 4.15, Prop. 4.233]{G} implies
that either $G_1=\PSL_2(\F_q)$ or $G_1=\SL_2(\F_q)$. Applying
Lemma \ref{sl2}, we arrive to the desired contradiction.
\end{proof}

\section{Proof of Lemmas \ref{double} and \ref{sl2}}
\label{pl}
 We start with some auxiliary constructions related to
the permutation groups \cite{Klemm,Mortimer,Ivanov}.

 Let $B$ be a finite set consisting of $n \ge 5$ elements. We
write $\Perm(B)$ for the group of permutations of $B$. A choice of
ordering on $B$ gives rise to an isomorphism $\Perm(B) \cong \Sn$.
 Let us assume that $n$ is {\sl odd} and consider the permutation module $\F_2^{B}$: the
$\F_2$-vector space of all functions $\varphi:B \to \F_2$. The
space $\F_2^{B}$ carries a natural structure of $\Perm(B)$-module
and contains  the stable hyperplane $Q_{B}:=(\F_2^{B})^0$ of
functions $\varphi$ with $\sum_{\alpha\in B}\varphi(\alpha)=0$.
 Clearly,
$Q_{B}$ carries a natural structure of faithful $\Perm(B)$-module.
For each permutation group $H\subset \Perm(B)$ the corresponding
$H$-module is called the {\sl heart} of the permutation
representation of $H$ on $B$ over $\F_2$
\cite{Klemm,Mortimer,Ivanov}.

\begin{lem}
\label{Kl} $\End_H(Q_{B})=\F_2$ if $n$ is odd and $H$ acts
$2$-transitively on $B$.
\end{lem}

\begin{proof}
See Satz 4 in \cite{Klemm}.
\end{proof}

\begin{proof}[Proof of Lemma \ref{double}]

Suppose $f(x)\in K[x]$ is a polynomial of odd degree $n \ge 5$
without multiple roots and $X:=J(C_f)$ is the jacobian of $C=C_f:
y^2=f(x)$. It is well-known that $g:=\dim(X)=\frac{n-1}{2}$.
  It is also
well-known (see for instance Sect. 5 of \cite{ZarhinTexel}) that
the image of $\Gal(K) \to \Aut(X_2)$ is isomorphic to $\Gal(f)$.
More precisely, let $\RR\subset K_a$ be the $n$-element set of roots of $f$, let $K(\RR)$ be
 the splitting field of $f$ and $\Gal(f)=\Gal(K(\RR)/K)$ the Galois group of $f$,
viewed as a subgroup of of the group $\Perm(\RR)$ of all
permutations of $\RR$. We have $\Gal(f) \subset \Perm(\RR)$. It is
well-known (see for instance, Th. 5.1 on p. 478 of
\cite{ZarhinTexel}) that  $\Gal(K) \to \Aut(X_2)$ factors through
the canonical surjection $\Gal(K) \twoheadrightarrow
\Gal(K(\RR)/K)=\Gal(f)$ and the $\Gal(f)$-modules $X_2$ and
$Q_{\RR}$ are isomorphic. In particular,
$$\End_{\Gal(K)}(X_2)=
\End_{\Gal(f)}(X_2)=\End_{\Gal(f)}(Q_{\RR}).$$ Assuming that
$\Gal(f)$ acts doubly transitively on $\RR$ and applying Lemma
\ref{Kl}, we conclude that
$$\End_{\Gal(f)}(X_2)=\End_{\Gal(f)}(Q_{\RR})=\F_2.$$
\end{proof}

\begin{rem}
The assertion of Lemma \ref{double} is implicitly contained in the
proof of Prop. 3 in \cite{Mori2}.
\end{rem}

\begin{proof}[Proof of Lemma \ref{sl2}]
It is known \cite[corollary on p. 4]{J} that $\Q[\PSL_2(\F_q)]$ is
a direct product of matrix algebras (for all power primes $q$). Since
$\ker(\SL_2(\F_q) \twoheadrightarrow \PSL_2(\F_q))$ is the only
proper normal subgroup in $\SL_2(\F_q)$,  it suffices to deal only
with the group $\SL_2(\F_q)$ with $q=5, g=2$ or $q=7, g=3$ and
consider only direct summands of $\Q[\SL_2(\F_q)]$ that correspond
(in the sense of Lemma 24.7 on p. 124 of \cite{DA}) to {\sl
faithful} irreducible characters of degree $q-1$ with values in
$\Q$.

  Let $\chi$ be an {\sl
irreducible faithful irreducible} character of degree $q-1$ with
values in $\Q$. Then (in the notations of \cite[\S 38]{DA})
$\chi=\theta_j$ where $j$ is an integer with $1 \le j \le
\frac{q-1}{2}$. If $z$ is the only nontrivial central element of
$\SL_2(\F_q)$ then $\theta_j(z)=(-1)^j (q-1)$.
 The faithfulness of $\chi$ implies (thanks to Lemma 2.19 \cite{Isaacs}) that  $\theta_j(z)\ne q-1$, i.e. $j$ is
odd. Let $b \in \SL_2(\F_q)$ be an element of order $q$ and
$\sigma$ a primitive $q+1$th root of unity. Then \cite[p. 228]{DA}
$$\chi(b)=\theta_j(b)=-(\sigma^j+\sigma^{-j}).$$

Assume that $q=7$. Then either $j=1$ or $j=3$. Also $q+1=8$ and we
may choose
$$\sigma=\frac{1+\sqrt{-1}}{\sqrt{2}}.$$
Then if $j=1$ then $\chi(b)=-\sqrt{2}$ and if $j=3$ then
$\chi(b)=\sqrt{2}$. In both cases $\chi(b)$ does not lie in $\Q$.
It follows that $\Q[\SL_2(\F_7)]$ does not have direct summands
isomorphic to the matrix algebras  of size $3$ over quaternion
$\Q$-algebras (including $\H_p$).

Assume that $q=5$. Then $j=1$ and $\chi=\theta_1$.  Then $q+1=6$
and the multiplicative order $n$ of $\sigma^j$ equals $6=2\cdot
3$. Also $\sigma^{2j}=\sigma^2$ is a primitive cubic root of
unity. Let $D$ be the direct summand of $\Q[\SL_2(\F_5)]$ attached
to $\chi$. It follows from the case (c) of Theorem on p. 4 of
\cite{J} (see also \cite[Th. 6.1(ii)]{F} (with
$\epsilon=\delta=1$)) that $D$ is isomorphic to to the matrix
algebra $\M_2(\H)$ where $H$ is a quaternion $\Q$-algebra ramified
(exactly) at $\infty$ and $3$. (This means that $H \cong \H_3$ and
$D \cong \M_2(\H_3)$.) It follows that if $D$ is isomorphic to
$\M_2(\H_p)$ then $p=3$.
\end{proof}

\section{Not supersingularity}

\label{main3p} We keep all the notations and assumptions of
Theorem \ref{supernot}.

We write $T_2(X)$ for the $2$-adic Tate
module of $X$ and
$$\rho_{2,X}:\Gal(F) \to \Aut_{\Z_2}(T_2(X))$$
for
the corresponding $2$-adic representation. It is well-known that $T_2(X)$ is
a free $\Z_2$-module of rank $2\dim(X)=2g$ and
$$X_2=T_2(X)/2 T_2(X)$$
(the equality of Galois modules). Let us put
$$H=\rho_{2,X}(\Gal(F)) \subset \Aut_{\Z_2}(T_2(X)).$$
Clearly, the natural homomorphism
$$\bar{\rho}_{2,X}:\Gal(F) \to \Aut(X_2)$$
defining the Galois action on the points of order $2$ is the composition of
$\rho_{2,X}$ and (surjective) reduction map modulo $2$
$$\Aut_{\Z_2}(T_2(X)) \to \Aut(X_2).$$
This gives us a natural (continuous) {\sl surjection}
$$\pi:H \to \bar{\rho}_{2,X}(\Gal(F)) \cong G,$$
 whose kernel consists of elements of $1+2 \End_{\Z_2}(T_2(X))$.
The choice of polarization on $X$ gives rise to a non-degenerate
alternating bilinear form (Riemann form) \cite{MumfordAV}
$$e: V_{2}(X) \times V_2(X) \to \Q_2(1) \cong \Q_2.$$
Since $F$ contains all $2$-power roots of unity, $e$ is $\Gal(F)$-invariant
and therefore is $H$-invariant. In particular,
$$H \subset \Sp(V_2(X),e)\subset\SL(V_2(X)).$$
Here $\Sp(V_2(X),e)$ is the symplectic group attached to $e$. In
particular, the $H$-module $V_2(X)$ is symplectic.

There exists a finite Galois extension $L$ of $F$ such that all
endomorphisms of $X$ are defined over $L$. Clearly,
$\Gal(L)=\Gal(F_a/L)$ is an open normal subgroup of finite index
in $\Gal(F)$ and
$$H'=\rho_{2,X}(\Gal(L)) \subset \Aut_{\Z_2}(T_2(X))\subset \Aut_{\Q_2}(V_2(X)))$$
is an open normal subgroup of finite index in $H$.
 We write $\End^0(X)$ for
the $\Q$-algebra $\End(X)\otimes\Q$ of endomorphisms of $X$.

There exists a finite Galois extension $L$ of $F$ such that all
endomorphisms of $X$ are defined over $L$. We write $\End^0(X)$
for the $\Q$-algebra $\End(X)\otimes\Q$ of endomorphisms of $X$.

Since $X$ is supersingular,
$$\dim_{\Q}\End^0(X)=(2\dim(X))^2=(2g)^2.$$
Recall (\cite{MumfordAV}) that the natural map
$$\End^0(X)\otimes_{\Q}\Q_{2} \to \End_{\Q_{2}}V_{2}(X)$$
is an embedding. Dimension arguments imply that
$$\End^0(X)\otimes_{\Q}\Q_{2} = \End_{\Q_{2}}V_{2}(X).$$
Since all endomorphisms of $X$ are defined over $L$, the image
$$\rho_{2,X}(\Gal(L)) \subset \rho_{2,X}(\Gal(F)) \subset\Aut_{\Z_2}(T_2(X))
\subset\Aut_{\Q_2}(V_2(X))$$ commutes with $\End^0(X)$. This
implies that $\rho_{2,X}(\Gal(L))$ commutes with
$\End_{\Q_{2}}V_{2}(X)$ and therefore consists of scalars. Since
$$\rho_{2,X}(\Gal(L)) \subset \rho_{2,X}(\Gal(F)) \subset \SL(V_2(X)),$$
$\rho_{2,X}(\Gal(L))$ is a finite group. Since $\Gal(L)$ is a
subgroup of finite index in $\Gal(F)$, the group
$H=\rho_{2,X}(\Gal(F))$ is also finite. In particular, the kernel
of the reduction map modulo $2$
$$\Aut_{\Z_2}T_2(X) \supset H \to G \subset \Aut(X_2)$$
consists of periodic elements and, thanks to Minkowski-Serre Lemma
\cite{SZ}, $Z:=\ker(\pi:H \to G)$ has exponent $1$ or $2$. In
particular, $Z$ is commutative. Since
$$Z \subset H \subset \Sp(V_2(X))\cong \Sp(2g,\Q_{2}),$$
$Z$ is a $\F_2$-vector space of dimension $\le g$.

Let $G_1$ be a minimal subgroup of $H$ such that $\pi(G_1)=G$.
(Since $H$ is finite, such $G_1$ always exists.) Since $G$ is
perfect, $G_1$ is also perfect. (Otherwise, we may replace $G_1$
by smaller $[G_1,G_1]$.) Clearly,
$$Z_1:=\ker(\pi: G_1 \twoheadrightarrow G) \subset Z$$
is also a $\F_2$-vector space of dimension $\le g$. We have
$$Z_1 \subset G_1 \subset H \subset \Sp(V_2(X))\cong \Sp(2g,\Q_{2}).$$
In particular, the symplectic $G_1$-module is a lifting of the
$G_1$($\twoheadrightarrow G$)-module $X_2$.

I claim that the natural representation of $G_1$ in the
$2g$-dimensional $\Q_2$-vector space $V_2(X)$ is absolutely
irreducible. Indeed, let us put
$$E:=\End_{G_1}(V_2(X))\subset \End_{\Q_2}(V_2(X)).$$
Clearly,
$$O_E = E \bigcap \End_{\Z_2}(T_2(X))\subset \End_{\Z_2}(T_2(X))$$
is a $\Z_{2}$-algebra that is a free $\Z_{2}$-module, whose
$\Z_{2}$-rank coincides with $\dim_{\Q_{2}}(E)$. Notice that
$O_E$ is a {\sl pure} $\Z_{2}$-submodule in
$\End_{\Z_2}(T_2(X))$, i.e. the quotient $\End_{\Z_2}(T_2(X))/O_E$
is a torsion-free (finitely generated) $\Z_{2}$-module and
therefore a free $\Z_{2}$-module of finite rank. It follows that
the natural map
$$O_E/2 O_E \to \End_{\Z_2}(T_2(X))/2 
\End_{\Z_2}(T_2(X))=\End_{\F_2}(X_2)$$ is an embedding. Clearly, the
image of $O_E/2 O_E$ in $\End_{\F_2}(X_2)$ lies in $\End_{G}(X_2)$. Since
$\End_{G}(X_2)=\F_2$, we conclude that the rank of the free
$\Z_{2}$-module $O_E$ is $1$, i.e. $\dim_{\Q_{2}}(E)=1$. This
means that $E=\Q_{2}$, i.e. the $G_1$-module $V_2(X)$ is
absolutely simple.

Let $\chi: G_1 \to \Q_{2}$ be the character of the absolutely
irreducible faithful representation of $G_1$ in $V_2(X)$. Clearly,
$\chi$ is a faithful (absolutely) irreducible character of degree
$2g$. We need to prove that $\chi(G_1) \subset \Q$.

Let $F_1\subset F_a$ be the subfield of invariants of the subgroup
$$\{\sigma \in \Gal(F)\mid \rho_{2,X}(\sigma) \in G_1\} \subset
\Gal(F).$$ Clearly, $F_1$ is a finite separable algebraic
extension of $F$ and
$$G_1=\rho_{2,X}(\Gal(F_1)).$$
 Clearly, the image
$$\bar{\rho}_{2,X}(\Gal(F_1))\subset \Aut(X_2)$$
coincides with
$$\pi\rho_{2,X}(\Gal(F_1))=\pi(G_1)=\pi_1(G_1)=G \subset
\Aut(X_2).$$

Let $L_1$ be the finite Galois extension of $F_1$ attached to
$$\rho_{2,X}:\Gal(F_1))\to \Aut(T_2(X)).$$
Clearly, $\Gal(L_1/F_1)=G_1$. In addition, all $2$-power torsion
points of $X$ are defined over $L_1$. It follows (see
\cite{Silverberg}) that all the endomorphisms of $X$ are defined
over $L_1$. On the other hand, I claim that the ring
$\End_{F_1}(X)$ of $F_1$-endomorphisms of $X$ coincides with $\Z$.
Indeed, there is a natural embedding
$$\End_{F_1}(X) \otimes \Z/2\Z \hookrightarrow
\End_{\Gal(F_1})(X_2)=\F_2$$ that implies that the rank of the
free $\Z$-module $\End_{F_1}(X)$ does not exceed $1$ and therefore
equals $1$, i.e. $\End_{F_1}(X)=\Z$.

Since all the endomorphisms of $X$ are defined over $L_1$, there
is a natural homomorphism
$$\kappa: G_1=\Gal(L_1/F_1) \to \Aut(\End(X))$$
%,\quad u \mapsto ^{\kappa(\sigma)}u, \quad u \in \End(X), \sigma \in
%G_1$$ 
such that
$$\End_{F_1}(X)=\{u\in \End(X)\mid \kappa(\sigma)u=u
\quad \forall \sigma \in \Gal(L_1/F_1)=G_1\}$$ and
$$\sigma(ux)=(\kappa(\sigma)u)(\sigma(x))\quad \forall x\in X(L_1), u \in \End(X),
\sigma \in \Gal(L_1/F_1)=G_1.$$ 
Further we write $^{\kappa(\sigma)}u$ for $\kappa(\sigma)(u)$.
Since $\End_{F_1}(X)=\Z$, we conclude that
$$\Z=\{u\in \End(X)\mid \ ^{\kappa(\sigma)}u=u
\quad \forall \sigma \in \Gal(L_1/F_1)=G_1\}.$$
Since all $2$-power torsion points
of $X$ defined over $L_1$,
$$\sigma(ux)=\ ^{\kappa(\sigma)}u (\sigma(x))\quad \forall x\in T_2(X), u \in \End(X),
\sigma \in G_1.$$
Since $\Aut(\End(X))\subset \Aut(\End^0(X))$, one may view $\kappa$ as 
$$\kappa: G_1=\Gal(L_1/F_1) \to \Aut(\End^0(X)),
\quad u \mapsto\ ^{\kappa(\sigma)}u, \quad u \in \End^0(X), \sigma \in
G_1$$ and we have
$$\Q=\{u\in \End^0(X)\mid \ ^{\kappa(\sigma)}u=u
\quad \forall \sigma \in \Gal(L_1/F_1)=G_1\}$$ and
$$\sigma(ux)=\ ^{k(\sigma)}u(\sigma(x))\quad \forall x\in V_2(X), u \in \End^0(X),
\sigma \in G_1.$$
Recall that
$$\End^0(X)\subset
\End^0(X)\otimes_{\Q}\Q_2=\End_{\Q_2}(V_2(X))$$ and
$$G_1 \subset \GL(V_2(X))=(\End_{\Q_2}(V_2(X)))^*.$$
It follows that
$$\sigma u \sigma^{-1}=\ ^{k(\sigma)}u \quad \forall u \in
\End^0(X), \sigma \in G_1.$$
 By
Skolem-Noether theorem, every automorphism of the central simple
$\Q$-algebra $\End^0(X)\cong\M_g(\H_p)$ is an inner one. This
implies that for each $\sigma \in G_1$ there exists $w_{\sigma}\in
\End^0(X)^*$ such that
$$\sigma u \sigma^{-1}=\ ^{k(\sigma)}u =w_{\sigma} u w_{\sigma}^{-1}
\quad \forall u \in \End^0(X).$$ Since the center of $\End^0(X)$
is $\Q$, the choice of $w_{\sigma}$ is unique up to multiplication
by a non-zero rational number. This implies that $w_{\sigma} w_{\tau}$
equals $w_{\sigma\tau}$ times a non-zero rational number.

Let us put
$$c'_{\sigma}= \sigma w_{\sigma}^{-1} \in (\End_{\Q_2}(V_2(X)))^*.$$
Clearly, each $c'_{\sigma}$ commutes with $\End^0(X)$ and
therefore with $\End^0(X)\otimes_{\Q}\Q_2=\End_{\Q_2}(V_2(X))$. It
follows that all $c'_{\sigma}$ are scalars, i.e. lie in $\Q_2^*
\I$. (Here $\I$ is the identity map on $V_2(X)$.) Clearly, the
image
$$c_{\sigma} \in \Q_2^* \I/\Q^* \I\cong {\Q_2}^*/\Q^*$$
of $c'_{\sigma}$ in ${\Q_2}^*/\Q^*$ does not depend on the choice
of $w_{\sigma}$. It is also clear that the map
$$G_1 \to {\Q_2}^*/\Q^*, \sigma \mapsto c'_{\sigma}$$
is a group homomorphism. Since $G_1$ is perfect and
${\Q_2}^*/\Q^*$ is commutative, this homomorphism is trivial, i.e.
$c_{\sigma}=1$ for all $\sigma \in G_1$. This means that
$$c_{\sigma} \in \Q^*  \I \quad \forall \sigma \in G_1$$
and therefore
$$\sigma={(c'_{\sigma})}^{-1}w_{\sigma}\in \End^0(X)^* \quad \forall \sigma \in G_1.$$
Recall \cite{MumfordAV} that if one view an element $u \in
\End^0(X)$ as linear operator in $V_2(X)$ then the characteristic
polynomial $P_u(t)$ of $u$ has rational coefficients; in
particular, the trace of $u$ is a rational number. It follows that
$\chi(G_1) \subset \Q$.

 Let $M$ be the image of
$$\Q[G_1] \to \End^0(X).$$
Clearly, $M \otimes_{\Q}\Q_2$ coincides with the image of
$$\Q_2[G_1] \to \End^0(X)\otimes_{\Q}\Q_2=\End_{\Q_2}(V_2(X)).$$

Since the $G_1$-module $V_2(X)$ is absolutely simple,
$$\Q_2[G_1] \to \End_{\Q_2}(V_2(X))$$
is surjective. This implies that
$$\dim_{\Q}(M)=\dim_{\Q}(\End^0(X))$$
and therefore, $M=\End^0(X)$, i.e. $\Q[G_1] \to \End^0(X)$ is
surjective. The semisimplicity of $\Q[G_1]$ allows us to identify
$\End^0(X)$ with a direct summand of $\Q[G_1]$.

If $\ell$ is a prime number that does not divide order of $G_1$
then it is well-known that the group algebra $\Q_{\ell}[G_1]$ is a
direct product of matrix algebras over (commutative) fields. It
follows that $p$ divides order of $G_1$. Since $\#(G_1)$ equals
$\#G$ times a power of $2$ and $p$ is odd, we conclude that $p$
divides $\#G$. In particular, $G_1$ contains an element $u$ of
exact order $p$. Since
$$u \in G_1 \subset \End^0(X) \subset \End_{\Q_2}(V_2(X)),$$
 $P_u(t)$ is  a polynomial of degree $2g$ with rational
coefficients and one of its roots is a primitive $p$th root of
unity. It follows that $P_u(t)$ is divisible in $\Q[t]$ by the
$p$th cyclotomic polynomial $\Phi_p(t)=\frac{t^p-1}{t-1}$. Since
the degree of $\Phi_p$ is $p-1$, we conclude that the degree $2g$
of $P_u(t)$ is greater or equal than $p-1$, i.e. $2g\ge p-1$.

Assume for a while that the $G$-module $X_2$ is very simple. Since
$G_1\to G$ is surjective, the $G_1$-module $X_2$ and its lifting
$V_2(X)$ are also very simple $G_1$-modules \cite[Remark 5.2
(i,v(a))]{ZarhinMMJ}. Since  $Z_1$ is normal in $G_1$, we
conclude, thanks to \cite[Remark 5.2(vii)]{ZarhinMMJ} that either
the $Z_1$-module $V_2(X)$ is absolutely simple or $Z_1$ consists
of scalars. Since $Z_1$ is a finite commutative group, it does not
admit absolutely irreducible representations of dimension $>1$.
Since $\dim_{\Q_2}(V_2(X))=2g>1$, we conclude that $Z_1$ consists
scalars; in particular, $Z_1$ is a central subgroup in $G_1$.
Since $$Z_1\subset G_1\subset \Sp(V_2(X))\cong \Sp(2g,\Q_{2}),$$
either $Z=\{1\}$ or $Z=\{\pm 1\}$. This implies that $Z_1$ is a
cyclic group of order $1$ or $2$.

Further we no longer assume that the $G$-module $X_2$ is very
simple. Assume instead that every homomorphism from $Z$ to
$\GL(g-1,\F_2)$ is trivial. I claim that in this case $Z$ is again
a central subgroup of $G_1$. Indeed, the short exact sequence
$$1 \to Z \hookrightarrow G_1 \twoheadrightarrow G \to 1$$
defines, in light of commutativeness of $Z$, a  natural
homomorphism
$$\eta:G \to \Aut(Z)$$
which is trivial if and only if $Z$ is central in $G_1$. Clearly, $\eta(G)$ is a finite perfect group.
Recall that $Z$ is an elementary $2$-group, i.e. $Z \cong \F_2^r$ for some nonnegative integer $r$.
Clearly, we may assume that $r \ge 1$ and therefore $\Aut(Z) \cong \GL(r,\F_2)$.
If $r \le g-1$ then we are done. Suppose that $r=g$. Then $Z$ must contain
$$\{\pm 1\} \subset \Sp(V_2(X).$$
Since $\{\pm 1\}$ is a central subgroup of $G_1$, the elements of
$\eta(G)\subset \Aut(Z)$ act trivially on $\{\pm 1\}$. Since the
quotient $Z/\{\pm 1\}$ has $\F_2$-dimension $g-1$, elements of
$\eta(G)$ act trivially on $Z/\{\pm 1\}$. This implies that
$\eta(G)$ is isomorphic to a subgroup of the commutative group
$\Hom(Z/\{\pm 1\}, \{\pm 1\}$. Since $\eta(G)$ is perfect, we
conclude that $\eta(G)=\{1\}$, i.e. $Z$ is a central subgroup and
therefore is either $\{1\}$ or $\{\pm 1\}$.

\section{Hyperelliptic two-dimensional jacobians in characteristic $3$}
\label{g2c3}
 Throughout this section $K$ is a field of
characteristic $p=3$ and $K_a$ its algebraic closure, $n=5$ or
$6$,
$$f(x)=\sum_{i=0}^n a_i x^i \in K[x]$$
a separable polynomial of degree $n$, i.e. all $a_i\in K, a_n \ne
0$ and $f$ has no multiple roots. We write $\Gal(f)\subset \SS_n$
for the Galois group of $f$ over $K$.

 Let $C_f$ be the
hyperelliptic curve $y^2=f(x)$ over $K_a$.

\begin{lem}
\label{ssS} Suppose that $n=\deg(f)=5$ and $a_4=0$.
\begin{itemize}
\item[(i)]
The jacobian $J(C_f)$ of $C_f$ is a supersingular abelian variety
over $K_a$ if and only if $a_1=a_2=0$, i.e.
$$f(x)=a_5 x^5+a_3 x^3+a_0.$$
If this is the case then $J(C_f)$ is isogenous but not isomorphic
to a self-product of a supersingular elliptic curve.
\item[(ii)]
Suppose that $a_0 \ne 0$ (e.g., $f(x)$ is irreducible over $K$)
and $J(C_f)$ is a supersingular abelian variety. Then
$\Gal(f)\subset \A_5$ if and only if $-1$ is a square in $K$, i.e. $K$ contains 
$\F_9$.
\end{itemize}
\end{lem}

\begin{proof}
Since $p=3$, $f(x)^{(p-1)/2}=f(x)$. Let us consider  the matrices
$$M:=\begin{pmatrix}a_{p-1}&a_{p-2}\\a_{2p-1}&a_{2p-2}\end{pmatrix}=
\begin{pmatrix}a_2&a_1\\a_5&0\end{pmatrix}, M^{(3)}:=
\begin{pmatrix}a_2^3&a_1^3\\a_5^3&0\end{pmatrix}.$$
 Extracting cubic roots
from all entries of $M$ one gets the Hasse--Witt/Cartier--Manin
matrix $M^{(3)}$ of $C$ (with respect to the standard basis in the
space of differentials of the first kind) \cite{Manin},
\cite{Yui},~\cite[p. 129]{Oort}. Recall (~\cite[p. 78]{Manin},
\cite{Nygaard},
 ~\cite[Th. 3.1]{Yui},~\cite[Lemma 1.1]{Oort}) that the jacobian
 $J(C)$ is a supersingular abelian surface not isomorphic to a product of two
 supersingular elliptic curves if and only if $M \ne 0$ but
$$M^{(3)} M =0.$$
Clearly, $M \ne 0$, because $a_5 \ne 0$. It is also clear that
$$\det(M^{(3)} M)=\det(M^{(3)})\det(M)=(-a_1^3 a_5^3)(-a_1
a_5)=a_1^4 a_5^4.$$ Hence, if $M^{(3)} M =0$ then $a_1=0$.

Suppose that $a_1=0$. Then
$$M=\begin{pmatrix}a_{2}&0\\a_{5}&0\end{pmatrix},
M^{(3)}=\begin{pmatrix}a_{2}^3&0\\a_{5}^3&0\end{pmatrix},
  M^{(3)} M=
\begin{pmatrix}a_{2}^4&0\\a_{5}^3 a_2&0\end{pmatrix}.$$
We conclude that $M^{(3)} M=0$ if and only if $a_1=a_2=0$. It
follows that $J(C)$ is a supersingular abelian surface if and only
if $a_1=a_2=0$. Since $M \ne 0$, the jacobian $J(C)$ is not
isomorphic to a product of two
 supersingular elliptic curves. This proves (i).

 In order to prove (ii), let us assume that $J(C_f)$ is
 supersingular, i.e.,
 $$f(x)=a_5 X^5+a_3 x^3 +a_0.$$
 We know that $a_0\ne 0, a_5 \ne 0$.
 Let us put
 $$h(x):=a_5^{-1} f(x)=x^5+b_3 x^3+ b_0$$
 where $b_3=a_3/a_5, b_0=a_0/a_5$. Clearly, $b_0\ne 0$ and the Galois groups of
 $f(x)$ and $h(x)$ coincide. So, it suffices to check that
 $\Gal(h)\subset \A_5$ if and only if $-1$ is a square in $K$.

 The derivative $h'(x)$ of $h(x)$ is $5  x^4=- x^4$.
 Let $\alpha_1, \cdots , \alpha_5$ be the roots of $h$. Clearly,
 $$\prod_{i=1}^5\alpha_i= -b_0.$$
 It is well-known that the Galois group of $h$ lies in the
 alternating group if and only if its  discriminant
$$D=\prod_{i<j}(\alpha_i-\alpha_j)^2$$
 is a square in $K$. On the other hand, it is also well-known that
$$\prod_{i=1}^5 h'(\alpha_i)=:R(h,h')=(-1)^{\deg(h)(\deg(h)-1)/2}D.$$
(Here $R(h,h')$ is the resultant  of $h$ and $h'$.)
 It follows that
$$R(h,h')= \prod_{i=1}^5(-
\alpha_i^4)=-(\prod_{i=1}^5\alpha_i)^4=- (-b_0)^4=- b_0^4$$ and
therefore $ D=- b_0^4$. Clearly, $D$ is a square in $K$ if and
only if $-1$ is a square in $K$.
\end{proof}

\begin{ex}[Counterexamples for $\A_5$ and $\SS_5$]
Let $k$ be an algebraically closed field of characteristic $p=3$.
Let $K=k(z)$ be the field of rational functions in variable $z$
with constant field $k$. We write $\overline{k(z)}$ for an
algebraic closure of $k(z)$. According to Abhyankar \cite{AbB},
the Galois group of the polynomial
$$h(x)=x^5-z x^2+1 \in k(z)[x]=K[x]$$
is $\A_5$ (see also ~\cite[\S 3.3]{SerreB}). It follows that the
Galois group of the polynomial
$$f(x)=x^5 h\left(\frac{1}{x}\right)=x^5-z x^3+1=\sum_{i=1}^5 a_i x^i$$
 is also $\A_5$. (Here
$a_5=1, a_4=a_2=a_1=0, a_3=-z, a_0=1$.)

 Let us consider the hyperelliptic curve
$$C: y^2=x^5-z x^3+1$$
of genus $2$ over $\overline{k(z)}$.
 It follows from  Lemma \ref{ssS} that the jacobian
 $J(C)$ of $C$ is a supersingular abelian surface that is {\sl not} isomorphic
 to a product of two supersingular elliptic curves. Hence
 $\End(J(C))$ is isomorphic to a certain order in the matrix algebra of size $2$
 over the quaternion $\Q$-algebra ramified exactly at $3$ and
 $\infty$. See ~\cite[proposition 2.19]{Oort}) for an explicit
 description of this order.

Assume now that $k$ is an algebraic closure of $\F_3$. Let us put
$$K_0=\F_3(z)\subset K=k(z)\subset \overline{k(z)}.$$
Clearly, $-1$ is {\sl not} a square in $K_0$ and $\overline{k(z)}$
is an algebraic closure of $K_0$.  Also, $f(x)\in K_0[x]$. An
elementary calculation (as in the proof of  Lemma \ref{ssS}(ii))
shows that the discriminant of $f(x)$ is $-1$. This implies
 that the Galois group of  $f(x)$ over $K_0$ does not lie in
$\A_5$. It follows that the Galois group of $f(x)=x^5-z x^3+1$
over $K_0$ is $\SS_5$. However, as we have already seen, the
jacobian of $y^2=x^5-z x^3+1$ is supersingular.
\end{ex}

\begin{thm}
\label{char3} Let $K$ be a field with $\fchar(K)=3$,
 $K_a$ its algebraic closure,
$f(x) \in K[x]$ an irreducible separable polynomial of  degree
$n=5$ or $6$. Let us assume  that the Galois group $\Gal(f)$ of
$f$ is the full symmetric group $\Sn$. Assume, in addition, that
$-1$ is a square in $K$, i.e. $K$ contains $\F_9$.

Let $C=C_f$ be the hyperelliptic curve $y^2=f(x)$. Let $J(C_f)$ be
its jacobian, $\End(J(C_f))$ the ring of $K_a$-endomorphisms of
$J(C_f)$. Then $\End(J(C_f))=\Z$.
\end{thm}

\begin{proof}
[Proof of Theorem \ref{char3}] Thanks to Remark \ref{odd}, we may
and will assume that $n=5$. We have
$$f(x)=\sum_{i=0}^5 a_i x^i \in K[x]$$
where all the coefficients $a_i \in K$ and $a_0 \ne 0$. Let us put
$$\gamma:=\frac{a_4}{5 a_0}, \quad h(x):=f(x-\gamma).$$
Clearly, $h(x) \in K[x]$ is an irreducible polynomial of degree
$5$ and $\Gal(h)=\Gal(f)=\SS_5$. It is also clear that if
$$h(x)=\sum_{i=0}^5 b_i x^i \in K[x]$$
then $b_4=0, \ b_5=a_5 \ne 0$. The substitution $x_1=x+\gamma, \
y_1=y$ establishes a $K$-birational isomorphism between
hyperelliptic curves $C=C_f:y^2=f(x)$ and $C_1=C_h:y_1^2=h(x_1)$
and induces an isomorphism of the jacobians $J(C_f)$ and $J(C_h)$.

Suppose that $\End(J(C_f))\ne \Z$. Then it follows from Theorem
2.1 of \cite{ZarhinMRL} that $J(C_f)$ is a supersingular abelian
variety. It follows that $J(C_h)\cong J(C_f)$ is also a
supersingular abelian variety. Applying Lemma \ref{ssS}(ii) to
$h$, we conclude that $\Gal(h)\subset \A_5$, because $-1$ is a
square in $K$. However, $\Gal(h)=\SS_5$. We obtained the desired
contradiction.
\end{proof}

\begin{ex}
Let $k$ be an algebraically closed field of characteristic $3$.
Let $K=k(z)$ be the field of rational functions in variable $z$
with constant field $k$. We write $\overline{k(z)}$ for an
algebraic closure of $k(z)$. Let $h(x) \in k[x]$ be a {\sl Morse
polynomial} of degree $5$. This means that the derivative $h'(x)$
of $h(x)$ has $\deg(h)-1=4$ distinct roots $\beta_1, \cdots
\beta_{4}$ and $h(\beta_i) \ne h(\beta_j)$ while $i\ne j$. (For
example, $x^5-x$ is a Morse polynomial.)  Then a theorem of
Hilbert (\cite[theorem~4.4.5, p.~41]{Serre}) asserts that the
Galois group of $h(x)-z$ over $k(z)$ is $\Sn$. Let us consider the
hyperelliptic curve
$$C: y^2=h(x)$$
of genus $2$ over $\overline{k(z)}$ and its jacobian $J(C)$. It
follows from Theorem \ref{char3} that $\End(J(C_f))=\Z$. (The case
of $h(x)=x^5-x$ was earlier treated by Mori \cite{Mori1}.)
\end{ex}

\section{A corollary}
 Combining Theorems \ref{main} and \ref{char3} together
with Theorem 2.3 of \cite{ZarhinMMJ} and Theorem 2.1 of
\cite{ZarhinMRL}, we obtain the following statement.

\begin{thm}
\label{general} Let $K$ be a field with $\fchar(K) \ne 2$,
 $K_a$ its algebraic closure,
$f(x) \in K[x]$ an irreducible separable polynomial of degree $n
\ge 5$ such that the Galois group of $f$ is either $\Sn$ or $\An$.
If $\fchar(K)=3$ and $n \le 6$ then we additionally assume that 
 $\Gal(f)=\Sn$ and $K$ contains $\F_9$.

Let $C_f$ be the hyperelliptic curve $y^2=f(x)$. Let  $J(C_f)$ be
its jacobian, $\End(J(C_f))$ the ring of $K_a$-endomorphisms of
$J(C_f)$. Then  $\End(J(C_f))=\Z$.
\end{thm}

\end{document}